\documentclass[12pt]{article}

\usepackage{amsmath}
\usepackage{amssymb}
\usepackage{amsthm}
\usepackage{tikz}
\usepackage{lineno}
\usepackage{enumerate}

\newcommand{\cprod}{\mathbin{\square}}
\newcommand{\kne}[2]{\ensuremath K_{#1:#2}}
\newcommand{\hc}[1]{\ensuremath \frac{1}{2}Q_{#1}}

\newtheoremstyle{plainsl}%
	{\topsep}
	{\topsep}
	{\slshape} 
	{}
	{\normalfont\bfseries}
	{.}
	{ }
	{}

\swapnumbers

{\theoremstyle{plainsl}
\newtheorem{theorem}{Theorem}[section]
\newtheorem{lemma}[theorem]{Lemma}
\newtheorem{corollary}[theorem]{Corollary}}
{\theoremstyle{remark}

}

\newcommand\cref[1]{Corollary~\ref{cor:#1}}

\def\sqr#1#2{{\vbox{\hrule height.#2pt
    \hbox{\vrule width.#2pt height#1pt \kern#1pt
        \vrule width.#2pt}\hrule height.#2pt}}}
\def\eqed{\sqr53}
\def\qed{%
    \ifmmode\eqno\eqed
    \else\nobreak\ \hfill\eqed\medbreak\fi}

\begin{document}

\title{Cores of Geometric Graphs}
\author{
	Chris Godsil and Gordon F. Royle
	\footnote{This work was supported by Canada's NSERC and completed while 
	the second author was visiting the Department of Combinatorics \& Optimization, 
	University of Waterloo.}\\[6pt]
	\small Combinatorics and Optimization, University of Waterloo\\[-0.8ex]
	\small \texttt{cgodsil@uwaterloo.ca}\\[3pt]
	\small Mathematics \& Statistics, University of Western Australia\\[-0.8ex]
	\small \texttt{gordon@maths.uwa.edu.au}
	}

\date{April 2008}

\maketitle

\begin{abstract}
Cameron and Kazanidis have recently shown that rank-3 graphs are either cores or have
complete cores, and they asked whether this holds for all strongly regular graphs. We 
prove that this is true for the point graphs and line graphs of generalized quadrangles and that
when the number of points is sufficiently large, it is also true for the block graphs of 
Steiner systems and orthogonal arrays. 
\end{abstract}

\section{Introduction}

A graph \textsl{homomorphism} between two graphs $X$ and $Y$ is a function
$$
\varphi: V(X) \rightarrow V(Y)
$$
such that $\varphi(x)\varphi(y)$ is an edge of $Y$ whenever $xy$ is an edge of $X$ (see 
\cite{MR1468789, MR2089014} for background material on graph homomorphisms). Unlike a graph isomorphism (which it superficially resembles) a homomorphism can map several vertices of
$X$ onto a single vertex of $Y$.  Two graphs $X$ and $Y$ are \textsl{homomorphically equivalent} if
there is a homomorphism from $X$ to $Y$ and a homomorphism from $Y$ to $X$.

A homomorphism from $X$ into the complete graph $K_q$ is equivalent to 
a $q$-colouring of $X$ and we can view homomorphisms as generalizations of 
colourings where the complete graphs are replaced by other families of graphs. For
example, if we use the Kneser graphs $\kne{v}{r}$ as the target graphs, then we get 
the theory of \textsl{fractional colourings}. 

An \textsl{endomorphism} of a graph $X$ is a homomorphism from $X$ to itself, and the set of
all endomorphisms of $X$ is denoted ${\rm End}(X)$. The composition of two endomorphisms
is again an endomorphism and the identity mapping is an endomorphism and so ${\rm End}(X)$ is
a \textsl{monoid}. Clearly any automorphism of $X$ is an endomorphism and so 
$$
{\rm Aut}(X) \subseteq {\rm End}(X).
$$
A \textsl{proper} endomorphism is an element $\varphi \in {\rm End(X)} \backslash {\rm Aut}(X)$, and in this situation the image of $\varphi$ is a proper subgraph of $X$. One of the fundamental concepts
in the theory of graph homomorphisms is that of a \textsl{core} which is a graph with no proper
endomorphisms. Cores play an important role in graph homomorphisms because every
graph is homomorphically equivalent to a unique core and thus the cores are the canonical
representatives of homomorphism equivalence classes. We let $X^\bullet$ denote the unique
core homomorphically equivalent to a graph $X$ and call this \textsl{the core of} $X$. The 
core of $X$ is necessarily isomorphic to an induced subgraph of $X$.

If the core of a graph is a complete graph, say $K_q$, then the graph not only contains 
a $q$-clique but is also $q$-colourable and hence 
$$
\omega (X) = \chi (X) = q
$$
where $\omega(X)$ is the size of the largest clique and $\chi(X)$ is the chromatic number of $X$. The converse of this statement also holds, in that if $\omega(X) = \chi(X)$ then the core of $X$ is necessarily complete. 

Cameron \& Kazanidis \cite{camkaz} considered the cores of symmetric graphs, in particular rank-3 graphs and proved that if a rank-3 graph does not have a complete core, then it is a core itself.

Rank-3 graphs are necessarily strongly regular and Cameron \& Kazanidis (personal communication) tentatively conjectured that \textsl{all} strongly regular graphs are cores or have complete cores. A result of Neumaier  \cite{MR564298} (extending a result of Bose) is that for any fixed minimum eigenvalue $\tau$, all but finitely many strongly regular graphs with an eigenvalue
equal to $\tau$ are the point graphs of certain partial geometries. Therefore these \textsl{geometric graphs} are fundamental classes of strongly regular graphs that are not (necessarily) rank-3 graphs.

In this paper, we show that Cameron \& Kazanidis's conjecture holds for the point graphs of generalized quadrangles and, provided the number of points is sufficiently large, the graphs arising from $2$-$(v,k,1)$ designs and orthogonal arrays.

\section{Geometric Graphs}

A \textsl{partial geometry} $PG(s,t,\alpha)$ is a point-line incidence structure satisfying the following
conditions:

\begin{enumerate}[(1)]
\item Two distinct lines meet in at most one point, and two distinct points are joined by at most one line.
\item Each line contains $s+1$ points, and each point has $t+1$ lines through it.
\item Given a point $P$ and a line $\ell$ not containing $P$ there are exactly $\alpha$ lines
through $P$ meeting $\ell$.
\end{enumerate}

The \textsl{dual} of a partial geometry $PG(s,t,\alpha)$ is obtained by exchanging the roles of 
points and lines and is a partial geometry $PG(t,s,\alpha)$.  The \textsl{point graph} (or \textsl{collinearity graph}) of a partial geometry is the graph whose vertex set is the point set of the partial geometry and where two vertices are adjacent if and only if the corresponding points are collinear in the partial geometry. If every two points are collinear, then the point graph is complete, and otherwise it is an $(n,k,\lambda,\mu)$ strongly regular graph with the following parameters:
\begin{align*}
n &=  (st + \alpha) (s+1) / \alpha,\\
k &=  (t+1)s\\
\lambda &= (s-1) +t(\alpha - 1),\\
\mu &= (t+1) \alpha.
\end{align*}
A strongly regular graph with these parameters has eigenvalues $k$, $\theta$ and $\tau$ with
multiplicities $1$, $m_\theta$ and $m_\tau$ respectively where
\begin{align*}
\theta &= s-\alpha, \\
m_\theta &= \frac{st(s+1)(t+1)}{\alpha (s+t+1-\alpha)},\\
\tau &= -1-t,\\
m_\tau &= \frac{(st+\alpha)(s+1-\alpha)}{\alpha(s+t+1-\alpha)}.
\end{align*}

Three of the most important  families of partial geometries are those arising from designs, orthogonal arrays and generalized quadrangles. They can be described as follows:

\begin{description} 

\item[Designs]

A $2$-$(v,k,1)$ design consists of a set of $v$ points together with a set of $k$-subsets of the points such that any two points are contained in a unique block. If we construct a geometry with points being the \textsl{blocks} of the design and the lines being the \textsl{points} of the design where incidence
is natural, then we get a partial geometry $PG(s,t,\alpha)$ with
$$
s = (v-k)/(k-1) \qquad t = k-1 \qquad \alpha = k = t+1.
$$

\item [Orthogonal arrays ]

An orthogonal array $OA(k,n)$ is a $k \times n^2$ array with entries in $\{1,\ldots,n\}$ such that 
no $2 \times n^2$ sub-array contains repeated columns. Given an $OA(k,n)$ we can
construct a geometry with points being the $n^2$ columns and lines being the
$kn$  ``row, symbol'' pairs $(r,s)$ where line $(r,s)$ is incident with the points that 
have symbol $s$ in row $r$. This is a partial geometry $PG(s,t,\alpha)$ with 
$$
s = n-1 \qquad t = k-1 \qquad \alpha = k-1 = t.
$$
The dual of this partial geometry is called a \textsl{transversal design}.

\item [Generalized quadrangles]

A generalized quadrangle  $GQ(s,t)$ is defined to be a partial geometry with $\alpha = 1$; the \textsl{classical} generalized quadrangles arise from the theory of polar spaces, but there
are many other classes known and a substantial literature devoted to them \cite{MR767454}.

\end{description}

We add a few explanatory remarks regarding this list:

\begin{enumerate}[(1)] 
\item It is clear that $1 \leq \alpha \leq \min(s+1,t+1)$ and so the three classes above have 
the maximum possible, second maximum possible and smallest possible values of
$\alpha$ respectively.  Between these extremes, the partial 
geometries with $1 < \alpha < \min(s,t)$ are known as \textsl{proper} partial geometries and there are a few infinite families and a number of sporadic examples known (Thas \cite{thas_handbook}).

\item For any given value of $k$, there are $2$-$(v,k,1)$ designs with arbitrarily large $v$ and 
$OA(n,k)$ with arbitrarily large $n$, thereby providing families of strongly regular graphs with 
arbitrarily many vertices and fixed minimum eigenvalue $-k$. Neumaier's result shows that these
two families are the \textsl{only} infinite families of strongly regular graphs with minimum eigenvalue $-k$.

\item The point graphs of the \textsl{duals} of the partial geometries defined above are complete 
for the designs and complete multipartite for the orthogonal arrays, hence we need not consider
them any further.
\end{enumerate}

\section{Cores of geometric graphs}

The point graph of a partial geometry ${\mathcal S}$ contains $(s+1)$-cliques consisting of the points on
a line of ${\mathcal S}$ and these are the largest possible cliques (unless the point graph is complete). If the partial geometry is ``large enough'' (in a sense to be made precise below) then there are
no other $(s+1)$-cliques.  As any endomorphism of the point graph must necessarily map
$(s+1)$-cliques onto $(s+1)$-cliques, we can view the endomorphism as mapping lines
onto lines in the partial geometry. We deal separately with the cases where $\alpha$ is 
large and small.

\begin{theorem}
\label{thm:bigalpha}
Let $X$ be the point graph of a partial geometry ${\mathcal S} = PG(s,t,\alpha)$ 
with $s, t > 1$ and $\alpha > (t+1)/2$ and suppose that all $(s+1)$-cliques of $X$ are 
lines of ${\mathcal S}$. 
Then every proper endomorphism of $X$ has a single $(s+1)$-clique as its image, and is
therefore a colouring of $X$.
\end{theorem}

\proof
Suppose that $\varphi$ is a proper endomorphism of $X$ and that ${\mathcal S}$ contains
two lines $\ell_1$ and $\ell_2$ such that $\varphi(\ell_1) = \varphi(\ell_2) = \ell$. Then any
line meeting both $\ell_1$ and $\ell_2$ in distinct points is also mapped to $\ell$ by $\varphi$. Our
aim is to show that \textsl{every} point of ${\mathcal S}$ is mapped onto $\ell$ by $\varphi$.

Suppose first that $\ell_1$ and $\ell_2$ are parallel (note that this cannot occur if $\alpha = t+1$), and
let $Q$ be a point not on $\ell_1$ or $\ell_2$. There are $\alpha$ lines through $Q$ meeting $\ell_1$ and another $\alpha$ meeting $\ell_2$ and so if $2\alpha > t+1$ there is a line $m$ through $Q$ meeting both $\ell_1$ and $\ell_2$. Therefore $\varphi(m) = \ell$ and as every point of $S$ is
either on $\ell_1$, $\ell_2$ or on neither, the result follows. 

Now suppose that $\ell_1$ and $\ell_2$ meet in a point $P$ and, as previously, let $Q$ be a point
not on $\ell_1$ or $\ell_2$. If $\alpha = t+1$ (the ``design case'') then there is a line through $Q$
meeting $\ell_1$ and $\ell_2$ in distinct points and the result follows. If $\alpha < t+1$ then the
condition on $\alpha$ is only sufficient to guarantee a line through $Q$ meeting $\ell_1$
and $\ell_2$ in distinct points for those points $Q$  not collinear with $P$, and so we
can only conclude that the lines not through $P$ are all mapped onto $\ell$. However in this 
situation there are necessarily two parallel lines not through $P$ and so by using these
two lines as $\ell_1$ and $\ell_2$ we revert to the previous case.\qed

Next we deal with the generalized quadrangles. In a generalized 
quadrangle, three pairwise collinear points are necessarily collinear and so it
is immediate that every $(s+1)$-clique in the point graph is a line of the generalized quadrangle.

\begin{theorem}
\label{thm:smallalpha}
Let $X$ be the point graph of a generalized quadrangle ${\mathcal S} = GQ(s,t)$ with $s, t > 1$. 
Then every proper endomorphism of $X$ has a single $(s+1)$-clique as its image, and is
therefore a colouring of $X$.
\end{theorem}

\proof
Suppose that $\varphi$ is a proper endomorphism of $X$ and that ${\mathcal S}$ contains
two lines $\ell_1$ and $\ell_2$ such that $\varphi(\ell_1) = \varphi(\ell_2) = \ell$. 
If $\ell_1$ and
$\ell_2$ are parallel then there is necessarily a line meeting both $\ell_1$ and $\ell_2$ that
is also mapped on to $\ell$, and so by replacing $\ell_2$ if necessary, we assume that
$\ell_1$ and $\ell_2$ meet in a point $P$.

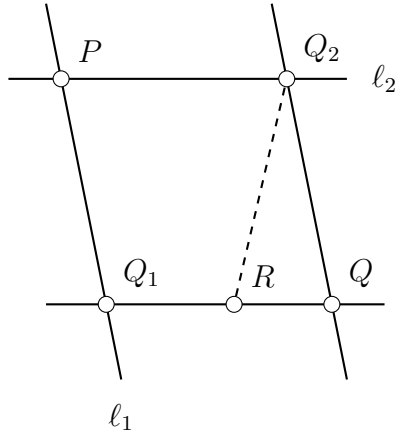
\begin{figure}
\begin{center}
\begin{tikzpicture}
\tikzstyle{geom}=[circle, draw=black, fill=white, inner sep=0.075cm]

\draw [thick] (0.5,1) -- (5,1);
\draw [thick] (0,4) -- (4.5,4);
\draw [thick] (0.5,5) -- (1.5,0);
\draw [thick] (3.5,5) -- (4.5,0);

\node [geom, label=45:$Q_1$] (Q1) at (intersection of 0.5,1--5,1 and 0.5,5--1.5,0) {};
\node [geom, label=45:${Q}$] (Q) at (intersection of 0.5,1--5,1 and 3.5,5--4.5,0) {};
\node [geom, label=45:${Q_2}$] (Q2) at (intersection of 0.5,4--5,4 and 3.5,5--4.5,0) {};
\node [geom, label=45:${P}$] (P) at (intersection of 0.5,4--5,4 and 0.5,5--1.5,0) {};
\node [geom, label=45:$R$] (R) at (3,1) {};
\node at (1.5,-0.5) {${\ell_1}$};
\node at (5,4) {${\ell_2}$};
\draw [thick, dashed] (R) -- (Q2);

\end{tikzpicture}
\end{center}
\caption{Impossible configuration in a generalized quadrangle}
\label{gqconfig}
\end{figure}

Now let $Q$ be a point not collinear with $P$ such that $\varphi(Q)$ is not 
on $\ell$.  Let $Q_1$, $Q_2$ be the unique points on 
$\ell_1$, $\ell_2$ respectively that are collinear with $Q$.  As $\varphi(Q)$ is not on 
$\ell$ it follows that $\varphi(Q_1) = \varphi(Q_2)$, because $X$ cannot contain
an $(s+1)$-clique with a point off the clique that is adjacent to more than one point of the clique. Now
consider another point $R$ on the line between $Q$ and $Q_1$ and consider the unique point $R_2$ on $\ell_2$ collinear with $R$. As $\varphi(Q)$ is not on $\ell$, nor is $\varphi(R)$ and so 
the argument above shows that $\varphi(Q_1) = \varphi(R_2)$. As an endomorphism is 
injective on an $(s+1)$-clique, it follows that $R_2 = Q_2$, which is not possible in a generalized
quadrangle (see Figure~\ref{gqconfig}).

Therefore every point not collinear with $P$ is mapped to a point on $\ell$ by $\varphi$. 
But any point
that \textsl{is} collinear with $P$ lies on a line containing at least two points 
\textsl{not} collinear
with $P$ and so this entire line is mapped on to $\ell$ by $\varphi$.\qed

Next we consider under what circumstances all the $(s+1)$-cliques of the line graph of a partial geometry ${\mathcal S} = PG(s,t,\alpha)$ are lines of ${\mathcal S}$.

\begin{theorem}
\label{thm:secondclique}
Let ${\mathcal S}$ be a partial geometry $PG(s,t,\alpha)$ with $\alpha > 1$ and let $X$ be its point graph. If $C$ is a clique in $X$ that is not contained in a line of ${\mathcal S}$ then
$$
|C| \leq 1 + (t+1)(\alpha-1).
$$
If equality holds then every line of the partial geometry contains either $0$ or $\alpha$ points
of $C$, and the points of $C$ together with the line intersections form a $2$-$(|C|,\alpha,1)$
design.
\end{theorem}

\proof
Let $P$ be a point in the clique $C$. There are $t+1$ lines through $P$ and every point of
$C$ lies on one of these lines. Each of these $t+1$ lines contains $P$ and at most $\alpha-1$
other points in $C$ (because a point off a line is collinear with only $\alpha$ points in total
on that line), and therefore $|C| \leq 1+ (t+1)(\alpha-1)$.

If equality holds then \textsl{every} line through $P$ contains $\alpha$ points of $C$, and as $P$ was
arbitrarily chosen this means that every line of ${\mathcal S}$ that is not disjoint from $C$ intersects
it in $\alpha$ points. It is immediate that the points of $C$ and the non-empty line intersections
form a $2$-design with the stated parameters.\qed

It is therefore clear that the hypotheses of Theorem~\ref{thm:bigalpha} are satisfied provided 
$s$ is sufficiently large.

\begin{corollary}
\label{cor:cliques_are_lines}
Let ${\mathcal S}$ be a partial geometry $PG(s,t,\alpha)$ and let $X$ be its point graph. If 
$$
s > (t+1)(\alpha-1)
$$
then every $(s+1)$-clique of $X$ is a line of ${\mathcal S}$.\qed
\end{corollary}

If $X$ is the block intersection graph of a $2$-$(v,k,1)$ design then the maximum size of
a non-line clique is $(k-1)^2 + (k-1) + 1$ which is the size of a projective plane of order $k-1$.  This bound is met if and only if the original design has a projective plane (of the same block size) as a subdesign. These cliques are strictly smaller than the lines whenever
$$
v > k(k^2-2k+2).
$$
For example, if we consider Steiner triple systems (i.e. $k=3$) then the line cliques are strictly
larger than all other cliques when $v > 15$. However there are several $2$-$(15,3,1)$ designs
that  contain Fano plane sub-designs.

If $X$ is the point graph of an $OA(k,n)$ then the maximum size of a non-line clique is
$(k-1)^2$ which is the size of an \textsl{affine} plane of order $(k-1)$. These non-line cliques
are strictly smaller than the lines whenever
$$
n > (k-1)^2.
$$
Again as an example, if we consider an $OA(3,n)$ (that is, a single Latin square of order $n$) 
then the line-cliques are strictly larger than all others whenever $n > 4$. For $n=4$ both of the
Latin squares have $2 \times 2$ sub-squares which yield $4$-cliques other than lines.

The two possibilities for the core allowed by Theorem~\ref{thm:bigalpha} and Theorem~\ref{thm:smallalpha} can both occur, and which of the two actually happens depends on the specific graph. For the graphs arising from designs, the core is complete if and only if the design is \textsl{resolvable}, while for the graphs arising from orthogonal arrays, the core is complete if and only if the $OA(k,n)$ can
be extended to an $OA(k+1,n)$. For generalized quadrangles the core is complete if and only 
if the point set can be partitioned into ovoids.

\section{More Cores}

In this section, we consider techniques for extending these results to a number of 
related classes
of graphs, including distance-transitive graphs and triangle-free distance regular graphs. 

In what follows we frequently use the fact that if $X$ is a graph and $\varphi$ is an endomorphism
onto its core, then we can assume that the image of $\varphi$ is a distinguished induced subgraph $Y$ isomorphic to $X^\bullet$ such that $\varphi$ restricted to $Y$ is the identity. Such an endomorphism is called a \textsl{retraction} from $X$ onto $Y$.

A graph is \textsl{distance transitive} if its automorphism group is transitive on ordered pairs of vertices at distance $i$ for all $i$. If $X$ is a distance-transitive graph, then there are constants $\{p_{ij}^k \mid 0 \leq i,j,k \leq d\}$ such that for any two vertices $v$ and $w$ at distance $k$, there are exactly $p_{ij}^k$ vertices at distance $i$ from $v$ and $j$ from $w$.

\begin{theorem}
If $X$ is a connected regular graph such that ${\rm Aut}(X)$ acts transitively on pairs of vertices at distance two, then  either $X$ is a core or it has a complete core.
\end{theorem}

\proof
Suppose that $\varphi$ is a proper retraction from $X$ onto its core $Y$. As $X$ is connected and regular and $\varphi$ is proper, there is some vertex $u$ such that $d_Y(\varphi(u)) < d_X(u)$ where $d_X(u)$ denotes the degree of $u$ in the graph $X$.
Therefore $u$ has two neighbours $v$ and $w$ such that $v$
and $w$ are not adjacent and $\varphi(v) = \varphi(w)$. Now suppose that $Y$ contains two
vertices at distance two, say $v'$ and $w'$. Then there is some automorphism $g \in {\rm Aut}(X)$ that maps$\{v',w'\}$ to $\{v,w\}$ and the composition $\varphi \circ g$ is a proper endomorphism of $Y$. As $Y$ has no proper endomorphisms this contradicts the choice of $v'$ and $w'$ and we conclude, that $Y$ is complete.\qed

\begin{corollary}
\label{cor:dtranscores}
A distance-transitive graph is either a core or has a complete core.\qed
\end{corollary}

The \textsl{Johnson graph} $J(v,k)$ has the $k$-subsets of a fixed $v$-set as its vertices, with two vertices adjacent if and only if the corresponding $k$-subsets meet in $k-1$ points. It is distance transitive and so by Corollary~\ref{cor:dtranscores}, it is either a core or has a complete core.  

\begin{theorem}
Let $X$ be a connected regular graph such that the neighbourhood of every vertex is either 
a thick generalized quadrangle $GQ(s,t)$ with $s, t > 1$ or the grid $K_p \cprod K_q$ with 
$1 < p \leq q$.  Then $X$ is either a core or has a complete core.
\end{theorem}

\proof
First suppose that the neighbourhood of each vertex is the grid $K_p \cprod K_q$, and that $\varphi$ is a retraction from $X$ onto its core $Y$. Then the valency in $Y$ of any vertex $v \in V(Y)$ is either $q$ or $pq$ depending on whether $\varphi$ restricted to $N(v)$ is a colouring or not. If some vertex $v \in V(Y)$ has valency $pq$ then its neighbourhood in $Y$ is $K_p \cprod K_q$ and so each vertex in \textsl{this} neighbourhood has valency in $Y$ at least $p+q-1$ and hence equal to $pq$. By connectivity, every vertex of $Y$ has valency $pq$ which contradicts the assumption that $\varphi$ is proper. We conclude that every vertex in $Y$ has degree $q$ and that $Y$ is complete. An analogous argument applies 
when the neighbourhoods are all thick generalized quadrangles.\qed 

As described above, it is often difficult to decide which of the two possibilities
holds for the core of a specific graph, although this can be done for certain
families of graphs. As an example, we consider the \textsl{halved cubes}, which are an
interesting family of graphs that arises naturally in the study of homomorphisms and colourings of Cayley graphs of elementary abelian 2-groups. 
The halved $n$-cube $\hc{n}$ has  the even weight vectors 
of ${\mathbb Z}_2^n$ as its vertex set, where two vertices are adjacent if and only if the corresponding vectors
differ in two places.  We need the following lemma, which may be found by combining two results in Hahn \& Tardif \cite {MR1468789}.

\begin{lemma}
\label{lem:transitivecores}
If $X$ is a vertex-transitive graph then its core $X^\bullet$ is vertex transitive and $|V(X^\bullet)|$ is
a divisor of $|V(X)|$. \qed
\end{lemma}

\begin{theorem}
If $n$ is a power of $2$ then the core of
$\hc{n}$ is the complete graph $K_n$ and otherwise $\hc{n}$ is a core.
\end{theorem}

\proof
By \cref{dtranscores}, either $\hc{n}$ is a core or its core is complete and therefore
a maximum clique of $\hc{n}$. As the maximum
cliques of $\hc{n}$ have $n$ vertices, Lemma~\ref{lem:transitivecores} shows that the core cannot be complete when $n$ is not a power of $2$. However, if $n$ is a power of 2, then the cosets of the extended Hamming code form an $n$-colouring of $\hc{n}$ and so $K_n$ is the core of $\hc{n}$.\qed

We now consider distance-regular graphs, which are graphs satisfying the same \textsl{combinatorial} regularity conditions as distance-transitive graphs, but which do not necessarily have automorphisms
underlying this regularity. More precisely, a graph of diameter $d$ is \textsl{distance regular} if there are constants $\{p_{ij}^k\mid 0 \leq i,j,k \leq d\}$ such that for any two vertices $v$, $w$ at distance $k$ there are exactly $p_{ij}^k$ vertices at distance $i$ from $v$ and $j$ from $k$. A distance-regular graph of diameter two is another term for strongly-regular graph, and the book by Brouwer, Cohen \& Neumaier \cite{MR1002568} is the definitive reference for larger diameters. 

We use a lemma that is proved in Godsil \& Royle \cite{MR1829620}; here a \textsl{$2$-arc} is a sequence $(u,v,w)$ of distinct vertices such that $u \sim v \sim w$.

\begin{lemma}
If $X$ is a connected non-bipartite graph such that every 2-arc lies in a shortest odd cycle, then $X$ is a core.\qed
\end{lemma}

\begin{corollary}
If $X$ is a distance-regular non-bipartite graph with no triangles, then $X$ is a core.
\end{corollary}

\proof
First observe that if $C$ is a shortest odd cycle of a graph $X$ then it is geodetic, i.e., the distance in $X$ between two vertices of $C$ is the distance between them in $C$. Now if $C$ has length $g > 3$, then $C$ contains two vertices $u$, $w$ at distance two such that there is a vertex at distance $(g-3)/2$ from $u$ and $(g-1)/2$ from $w$ in $C$.  Now if $(x,y,z)$ is an arbitrary $2$-arc of $X$ then $x$ and
$z$ are at distance two and so by distance regularity, there is at least one vertex $v$ at distance $(g-3)/2$ from $x$ and distance $(g-1)/2$ from $z$. The shortest paths from $v$ to $x$ and from $v$ to $z$ cannot have any vertex in common other than $v$ and cannot use $y$. Hence these two paths, together with the $2$-arc, form a shortest odd cycle.\qed

\section{Concluding Remarks}

We have examined all the strongly regular graphs with at most 36 vertices by computer 
and they all have the property that they have no proper endomorphisms other than colourings; although this may be viewed as supporting Cameron \& Kazinidis's conjecture, the number of strongly regular graphs increases dramatically with only a small increase in the number of vertices. 

A natural class of graphs to consider for possible counterexamples to the conjecture are the geometric graphs arising from partial geometries that \textsl{do not} satisfy Corollary~\ref{cor:cliques_are_lines}. In this case, the graphs may have ``extra'' $(s+1)$-cliques that do not arise from lines of the partial geometry and hence it may be possible to find an endomorphism whose image uses some of these extra cliques. There are several known $2$-$(15,3,1)$ designs with extra $7$-cliques (Fano subplanes) but none of these are counterexamples to the conjecture.

\end{document}